\documentclass{amsart}
\usepackage{amsmath}
  \usepackage{paralist}
  \usepackage{graphics} %% add this and next lines if pictures should be in esp format
  \usepackage{epsfig} %For pictures: screened artwork should be set up with an 85 or 100 line screen
\usepackage{graphicx}  \usepackage{epstopdf}%This is to transfer .eps figure to .pdf figure; please compile your paper using PDFLeTex or PDFTeXify.
 \usepackage[colorlinks=true]{hyperref}
\hypersetup{urlcolor=blue, citecolor=red}

\newtheorem{theorem}{Theorem}[section]

\theoremstyle{definition}

\usepackage{amssymb}
\newcommand{\e}{\varepsilon}
\newcommand{\R}{\mathbb{R}}

\newcommand{\N}{\mathbb{N}}
\newcommand{\Z}{\mathbb{Z}}
\newcommand{\T}{\mathbb{T}}
\DeclareMathOperator{\divv}{div}

\title[Non-uniqueness for the transport equation] %Use the shortened version of the full title
      {On some recent results concerning non-uniqueness for the transport equation}

\author[Stefano Modena]{}

 \email{Stefano.Modena@math.uni-leipzig.de}

\thanks{These notes were written during the visit of the author  to the Hausdorff Research Institute for Mathematics (HIM), University of Bonn, Jan-Apr 2019. This visit was supported by the HIM. Both, this support and the hospitality of HIM, are gratefully acknowledged.}

\begin{document}
\maketitle

% Enter the first author's name and address:
\centerline{\scshape Stefano Modena}
\medskip
{\footnotesize
% please put the address of the first author
 \centerline{Mathematisches Institut}
   \centerline{Universit\"at Leipzig}
   \centerline{D-04109 Leipzig, Germany}
} % Do not forget to end the {\footnotesize by the sign }

\bigskip

% The name of the associate editor will be entered by an editorial staff
% "Communicated by the associate editor name" is not needed for special issue.
% \centerline{(Communicated by the associate editor name)}

%The abstract of your paper
\begin{abstract}
In these notes we present some recent results concerning the non-uniqueness of solutions to the transport equation, obtained in collaboration with Gabriel Sattig and L\'aszl\'o Sz\'ekelyhidi in \cite{Modena2018, ModSze:2018, ModSat:2019}.
\end{abstract}

\section{Introduction}

These notes concern the problem of (non)uniqueness of solutions to the transport equation in the periodic setting
\begin{align}
\label{e:transport}
	\partial_t\rho+u\cdot \nabla\rho &=0,\\
\label{e:initial:cond}
	\rho_{|t=0}&=\rho^0
\end{align}
where $\rho:[0,T]\times \T^d\to\R$ is a scalar density, $u:[0,T]\times \T^d\to\R^d$ is a given vector field and $\T^d = \R^d / \Z^d$ is the $d$-dimensional flat torus. 

Unless otherwise specified, we assume in the following that $u \in L^1$ is \emph{incompressible}, i.e. 
\begin{equation}\label{e:incompressible}
	\divv u=0
\end{equation}
in the sense of distributions. Under this condition, \eqref{e:transport} is formally equivalent to the continuity equation
\begin{equation}\label{eq:continuity}
\partial_t \rho + \divv(\rho u) 	= 0.
\end{equation}

It is well known that the theory of classical solutions to \eqref{e:transport}-\eqref{e:initial:cond} is closely connected to the ordinary differential equation
\begin{equation}\label{e:ODE}
\begin{split}
	\partial_t X(t,x) &=u(t, X(t,x)),\\
	X(0,x) &= x.
\end{split}	
\end{equation}
More precisely, if $u$ is at least Lipschitz continuous, the solution to \eqref{e:transport}-\eqref{e:initial:cond} is given by the formula
\begin{equation}
\label{eq:ode:pde}
\rho(t, X(t,x))=\rho^0(x).
\end{equation}

There are several PDE models, related, for instance, to fluid dynamics or to the theory of conservation laws (see for instance \cite{DiPerna:1989:Annals,Crippa:2015er,LeBris2008,Lions:1996vo,Lions:1998vp}), where one has to deal with vector fields which are not  Lipschitz, but have lower regularity and therefore it is important to investigate the well-posedness of \eqref{e:transport}-\eqref{e:initial:cond} in the case of non-smooth vector fields.

There are several possibilities to state the well-posedness problem for \eqref{e:transport}-\eqref{e:initial:cond} in a weak setting; we describe now one possible way. Fix an exponent $p \in [1, \infty]$ and denote by $p'$ its dual H\"older, $1/p + 1/p' = 1$. The following two questions are of interest.
\begin{enumerate}[(a)]
\item Do existence and uniqueness of solutions to \eqref{e:transport}-\eqref{e:initial:cond} hold in the class of densities 
\begin{equation}
\label{eq:density:class}
\rho \in L^\infty(0,T; L^p(\T^d)) =: L^\infty_t L^p_x
\end{equation}
for a given vector field 
\begin{equation}
\label{eq:vf:class}
u \in L^1(0,T; L^{p'}(\T^d)) =: L^1_t L^{p'}_x?
\end{equation}
\item Is the relation \eqref{eq:ode:pde}, which links the PDE \eqref{e:transport} to the ODE \eqref{e:ODE} (or, in other words, the Eulerian world to the Lagrangian one) still valid, in some weak sense?
\end{enumerate}

Let us briefly comment on the choice of the classes \eqref{eq:density:class}-\eqref{eq:vf:class} for the density and the vector field, respectively. The choice of the class \eqref{eq:density:class} for the density is dictated by the following consideration. For smooth solutions to \eqref{e:transport}, every (spatial) $L^p$ norm remains constant in time. It is therefore natural in the weak setting to look for densities whose $L^p$ norm, if not constant, at least remains uniformly bounded in time. Once the class for $\rho$ is fixed, the choice \eqref{eq:vf:class} of the class for the vector field $u$ is as well natural, since in this way the product $\rho u \in L^1((0,T) \times \T^d)$ and hence the notion of distributional solution to \eqref{eq:continuity} (and thus also to \eqref{e:transport}) makes sense. 

\bigskip

This is the plan of these notes. In Section \ref{s:literature} we give a brief presentation of some well-posedness results and some counterexamples to well-posedness which can be found in the literature. In Section \ref{s:thm} we state the main theorem of these notes, Theorem \ref{thm:main}. In Section \ref{s:proof} we make some comments on the proof of Theorem \ref{thm:main}. 

We wish to stress that the aim of these notes is to give an informal presentation of some recent results concerning  non-uniqueness of solutions to the transport equation. For this reason, we intentionally avoid technicalities, we are quite vague in many points, many references are missing, and the statement of the main theorem is not presented in its full generality. For a more detailed discussion, we refer to \cite{Modena2018, ModSze:2018}.

\section{Well-posedness for the Cauchy problem in the weak setting}
\label{s:literature}

We sketch in this section a (far from complete) overview of the literature concerning the answers to questions (a) and (b) above. 

First of all, we remark that \emph{existence} of weak solutions in the class \eqref{eq:density:class}, for a given vector field as in \eqref{eq:vf:class}, is not a serious issue, because of the linearity of the transport equation. Indeed, to produce a weak solution to \eqref{e:transport}-\eqref{e:initial:cond}, it is enough to regularize the vector field $u$ and the initial datum $\rho^0$, to solve the regularized smooth problem and use the uniform bound in $L^\infty_t L^p_x$ to get a weakly converging sequence. By the linearity of the equation \eqref{e:transport}, the limit of such sequence is a weak solution to \eqref{e:transport}-\eqref{e:initial:cond}.

The big issue is thus \emph{uniqueness of weak solutions} and the \emph{relation \eqref{eq:ode:pde} between Eulerian and Lagrangian world}. 

\subsection{Uniqueness results}
\label{ss:uniqueness}

The first uniqueness result we mention is the celebrated theorem by DiPerna and Lions in 1989 \cite{DiPerna:1989vo}, when they proved that, if the vector field $u$, in addition to the integrability condition \eqref{eq:vf:class}, enjoys also the Sobolev regularity 
\begin{equation}
\label{eq:sobolev}
u \in L^1_t W^{1,p'}_x,
\end{equation}
then uniqueness of solutions holds in the class of densities \eqref{eq:density:class}. Let us remark that Di-Perna and Lions' Theorem is still true, even when the incompressibility condition \eqref{e:incompressible} is substituted by the weaker condition
\begin{equation}
\label{eq:div:bounded}
\divv u \in L^\infty_{tx}.
\end{equation}

Di-Perna and Lions' Theorem was extended in 2004 by Ambrosio \cite{Ambrosio:2004cva}, where he proved that, in the class of bounded densities (i.e. $p=\infty$ in our notation), uniqueness of solutions holds, if 
\begin{equation}
\label{eq:u:bv}
u \in L^1_t BV_x.
\end{equation}
Again, also Ambrosio's Theorem holds if \eqref{e:incompressible} is replaced by \eqref{eq:div:bounded}. 

Very recently, Bianchini and Bonicatto further extended Ambrosio's uniqueness result to vector fields which satisfy \eqref{eq:u:bv} and are \emph{nearly incompressible}. We do not want to enter into details here and to give a precise definition of \emph{near-incompressibility}. We only mention that such notion is the natural generalization of  \eqref{eq:div:bounded}, in the framework of $BV$ vector fields.  

We add two remarks to this list of results. The first one is the following. The proofs of the mentioned results are very subtle and involve several deep ideas and sophisticated techniques. We could however try to summarize the heuristics behind all of them as follows: (very) roughly speaking, a Sobolev or BV vector field $u$ is Lipschitz-like (i.e. $Du$ is bounded) on a large set and there is just a small ``bad'' set, where $Du$ is very large. On the big set where $u$ is ``Lipschitz-like'', the classical Cauchy-Lipschitz theory applies. Non-uniqueness phenomena could thus occur only on the small ``bad'' set. Uniqueness of solutions in the class of bounded densities (or $L^p$ densities, where $p$ is exactly the dual H\"older to the integrability exponent of $Du$, see \eqref{eq:sobolev}) is then a consequence of the fact that a \emph{bounded} (or $L^p$) density $\rho$ can not ``see'' this bad set, or, in other words, can not \emph{concentrate} on this bad set. 

A second interesting remark is that, roughly speaking, whenever uniqueness for the PDE \eqref{e:transport} holds \emph{in the class of bounded densities} (i.e. $p=\infty$) for a given vector field $u$, a uniqueness statement holds (in the sense of \emph{regular Lagrangian flow}, a notion we will not introduce in these notes, for a precise definition we refer, for instance, to \cite{Ambrosio2017}) also for the ODE \eqref{e:ODE} with the same vector field $u$. This can be seen, observing that the inverse flow map $\Phi(t) := X(t)^{-1} : \T^d \to \T^d$ is (at least in the smooth case) a bounded solution to \eqref{e:transport} with (vector valued) initial datum $\Phi(0,x) = x$.

\subsection{Non-uniqueness results}

From the analysis in the previous section it follows that the uniqueness results present in the literature concern vector fields
\begin{enumerate}[(a)]
\item which enjoy some form of exact or approximate \emph{incompressibility} (e.g. they have bounded divergence or they are nearly incompressible);
\item \emph{and} which are \emph{at least once differentiable} (in some weak sense, e.g. they are Sobolev or $BV$).
\end{enumerate}  

The counterexamples to uniqueness which can be found in the literature are, in general, based on the failure of at least one of these two conditions. For instance, already in the paper \cite{DiPerna:1989vo} by  Di-Perna Lions, it is possible to find an example of a Sobolev vector field with unbounded divergence and another example of an incompressible vector field which belongs to $L^1_t W^{1,s}_x$ for every $s<1$, but not to $L^1_t W^{1,1}_x$,  for which uniqueness of solutions fails. A further counterexample can be found in \cite{Depauw:2003wl} (an incompressible vector field which belongs to $L^1(\e, T; BV_x)$ for every $\e>0$ but not to $L^1(0,T; BV_x)$).

Let us also remark that the counterexamples mentioned so far are based on vector fields for which the associated ODE \eqref{e:ODE} has a degenerate behavior and therefore the Eulerian non-uniqueness is a consequence of the Lagrangian one.

\section{Statement of the main theorem}
\label{s:thm}

We mentioned in the previous section several uniqueness and non-uniqueness results and we observed that, in order to have uniqueness, the vector field $u$ must have some incompressiblity property and must possess one full spatial derivative. There is however one question we did not answer so far:

\smallskip
\begin{center}
for fixed $p \in [1, \infty)$, does uniqueness of solutions hold in the class of densities $L^\infty_t L^p_x$ for a given incompressible $u \in L^1_t W^{1, \tilde p}$, with $\tilde p < p'$? 
\end{center}
\smallskip

\noindent Recall that $p'$ is the dual H\"older exponent to $p$ and thus, if $\tilde p \geq p'$, then DiPerna-Lions' theory \cite{DiPerna:1989vo} guarantees uniqueness of solutions in $L^\infty_t L^p_x$. 

The answer to such question is not trivial at all. There are indeed two competing mechanisms, one playing for uniqueness, the other one playing against. 

On one side,  the incompressibility and the Sobolev regularity of $u$ imply uniqueness in the class of bounded densities (more precisely, in $L^1_t L^{\tilde p'}_x$, with $\tilde p'$ the dual H\"older to $\tilde p$) and thus, as observed in the previous section, uniqueness of solutions to the ODE \eqref{e:ODE} holds, in the sense of the regular Lagrangian flow. The Lagrangian picture is very well behaved.

On the other side, if ``$p$ is too small compared $\tilde p$'', it could happen (referring to the heuristics introduced in the previous section) that ``an $L^p$ density does see the \emph{bad set} of the $W^{1, \tilde p}$ vector field $u$'' and thus ``purely Eulerian'' non-uniqueness phenomena could occur. 

The following theorem, which is the main result we present in these notes, provides an answer to the question asked above.

\begin{theorem}[M., Sattig, Sz\'ekelyhidi]
\label{thm:main}
Let $p \in [1, \infty)$, $\tilde p \in [1, \infty)$. If 
\begin{equation}
\label{eq:condition}
\frac{1}{p} + \frac{1}{\tilde p} > 1 + \frac{1}{d},
\end{equation}
then there exist infinitely many incompressible vector fields 
\begin{equation*}
u \in C_t L^{p'}_x \cap C_t W^{1,\tilde p}_x
% \big( [0,T]; L^{p'}(\T^d; \R^d) \big) \cap C \big([0,T] ; W^{1, \tilde p} (\T^d; \R^d) \big), 
\end{equation*}
for which uniqueness of solutions to the transport equation \eqref{e:transport}
%\begin{equation*}
%\partial_t \rho + \nabla \rho \cdot u = 0, \quad \divv u = 0, % \quad \rho|_{t=0} = \rho^0, \quad \rho|_{t=1} = \rho^1
%\end{equation*}
fails in the class of densities $\rho \in C_t L^p_x$.  Moreover: 
\begin{itemize}
\item  if $p=1, p' = \infty$, then $u \in C \big( [0,T] \times \T^d \big) \cap C_t W^{1, \tilde p}_x$; 
\item the same result holds if the transport equation \eqref{e:transport} is replaced by the transport-diffusion equation
\begin{equation}
\label{eq:diffusion}
\partial_t \rho + \nabla \rho \cdot u = \Delta \rho %, \quad \divv u = 0,
\end{equation}
if, in addition, $p' < d$. 
\end{itemize}
\end{theorem}

\noindent Let us add some comments on the statement of Theorem \ref{thm:main}. 
\begin{enumerate}
\item The case $p = \infty$ is not considered. Indeed $p=\infty$ corresponds to the case of bounded densities and we have observed in Section \ref{ss:uniqueness} that, in this case, uniqueness holds even for $BV$ vector fields.
\item Similarly, also the  case $\tilde p = \infty$ is not considered. Indeed $\tilde p = \infty$ corresponds to the case of a Lipschtiz continuous vector field $u$ and, in this case, the classical Cauchy-Lipschitz theory for the ODE \eqref{e:ODE} provides a solution to \eqref{e:transport}-\eqref{e:initial:cond}, via the formula \eqref{eq:ode:pde}. 

\item In the case $p=1$, $p'=\infty$ (which correspond to $\tilde p < d$), the vector fields we construct are \emph{continuous}, not only \emph{bounded}. This shows that, in general, even the continuity of the vector field, in addition to the incompressibility and the Sobolev regularity, is not enough to guarantee uniqueness of weak solutions (compare with the result in \cite{Caravenna:2016kg, CarCri2018}). 

\item For the vector fields provided by Theorem \ref{thm:main}, uniqueness for the ODE \eqref{e:ODE} holds (in the sense of \emph{regular Lagrangian flow}): nevertheless, the PDE \eqref{e:transport} displays anomalous behavior. This shows that, for such vector fields, the relation between the Lagrangian and Eulerian world, summarized in Equation \eqref{eq:ode:pde}, is completely destroyed. This is even more evident in the case $p=1$, where the vector fields we construct are continuous and thus the trajectories of the regular Lagrangian flow are classical $C^1$ curves solving \eqref{e:ODE}.

\item In general, for the transport-diffusion equation \eqref{eq:diffusion} much stronger uniqueness results hold than for the transport equation \eqref{e:transport}. Indeed, the diffusion term $\Delta \rho$ is usually dominating (being the highest order term) and thus its regularizing effect translates, through the energy estimate, into a uniqueness statement for \eqref{eq:diffusion}. On the contrary, for the vector fields provided by Theorem \ref{thm:main}, the non-uniqueness generated by the first order term $\nabla \rho \cdot u$ is so strong that it beats even the second order term $\Delta \rho$. 

\end{enumerate}

\section{Some comments on the proof}
\label{s:proof}

We conclude these notes with some comments on the proof of Theorem \ref{thm:main}. 
Referring again to the heuristics introduced in Section \ref{ss:uniqueness}, the basic idea behind the proof of Theorem \ref{thm:main} is to ``concentrate the density $\rho$ on the \emph{bad} set of the vector field $u$''. 

This is done through a convex integration scheme, in the spirit of the papers by De Lellis, Sk\'ekelyhidi and collaborators on the Euler equations (see, in particular, \cite{DeLellis2014}). More precisely, the linear (in $\rho$) PDE \eqref{eq:continuity} is treated as a nonlinear PDE with both $\rho$ and $u$ as unknowns. The density $\rho$ and the field $u$ are constructed as limit of sequences
\begin{equation}
\label{eq:sequences}
\rho = \lim_{q \to \infty} \rho_q, \qquad u = \lim_q u_q, \qquad q \in \N,
\end{equation}
where the limits have to be taken in suitable norms and $(\rho_q, u_q)$ are approximate solutions to the transport equation, i.e.
\begin{equation}
\label{eq:cont:diff}
\partial_t \rho_q + \nabla \rho_q \cdot u_q  = \text{Error}_q, \qquad \divv u_q = 0,
\end{equation}
with $\text{Error}_q$ converging weakly to zero, as $q \to \infty$. 

The sequences $(\rho_q)_q, (u_q)_q$ are constructed recursively: assuming $\rho_q, u_q$ are given, as a first attempt, one defines
%\begin{equation*}
%\rho_{q+1} = \rho_q + a_q(t,x) \underbrace{{\Theta(\lambda_q x)}}_{\textit{fast  oscillation}}, \quad u_{q+1} = u_q + b_q(t,x) \underbrace{{W(\lambda_q x)}}_{\textit{fast oscillation}},
%\end{equation*} 
\begin{equation}
\label{eq:induction}
\rho_{q+1} = \rho_q + a_q(t,x) \Theta(\lambda_q x), \qquad u_{q+1} = u_q + b_q(t,x) W(\lambda_q x).
\end{equation} 
Here:
\begin{itemize}
\item $\lambda_q \in \N$ is an \emph{oscillation  parameter}, with $\lambda_q \to \infty$ as $q \to \infty$;
\item  $\Theta : \T^d \to \R$, $W: \T^d \to \R^d$ are fixed smooth profiles, called \emph{Mikado density} and \emph{Mikado field}, in the same spirit of the Mikado flows introduced by Daneri and Sz\'ekelyhidi in \cite{SzekelyhidiJr:2016tp} for the Euler equations; for a precise definition of $\Theta$ and $W$ we refer to the paper \cite{Modena2018};
\item $a_q, b_q$ are ``slow oscillating'' amplitudes, defined at each step in order to reduce $\text{Error}_q$ and to get, in the limit, a solution $(\rho, u)$ to \eqref{eq:continuity}. 
\end{itemize}

As in the framework of the Euler equations, the basic idea of convex integration is to choose the oscillation parameter $\lambda_q$ bigger and bigger along the iteration, and to use oscillations in order to reduce the error in \eqref{eq:cont:diff}.

The main difference between Theorem \ref{thm:main} and the theorems proven in the framework of the Euler equations (e.g. \cite{DeLellis2014, ise:18, del:ons:18}) is the following: in Theorem \ref{thm:main} we want to construct a vector field which is in $W^{1, \tilde p}_x$, i.e. it possesses one full derivative (in some $L^{\tilde p}$ space), whereas in the framework of the Onsager's conjecture for the Euler equations, the aim was to show the existence of anomalous $C^\gamma$ solutions, for every $\gamma < 1/3$, i.e. solutions which possess ``just $1/3$ of derivative'' (measured in a $\sup$ norm). 

How can we thus get such a $W^{1, \tilde p}$ bound? If a scheme as in \eqref{eq:induction} is used, one can easily see that problems arise. Indeed,  in order to have convergence of $D u_q$ in $L^{\tilde p}$, one should be able to provide a good bound of the distance $\|D u_{q+1} - Du_q\|_{L^{\tilde p}}$. However we have
\begin{equation}
\label{eq:estimate:der}
\|D u_{q+1} - Du_q\|_{L^{\tilde p}} \approx \lambda_q \|b_q\|_{L^\infty} \|DW\|_{L^{\tilde p}} 
\end{equation}
and the presence of the multiplicative factor $\lambda_q$ prevents the convergence of $D u_q$ in $L^{\tilde p}$. 

This issue can be solved using a \emph{concentration} argument, in the same spirit of what Buckmaster and Vicol did in the framework of the Navier-Stokes equations in their remarkable recent work \cite{Buckmaster:2017wf}, using \emph{intermittent Beltrami flows}.  In order to explain how the concentration argument works, let us think, for the time being, to the fixed Mikado density $\Theta$ and field $W$ as compactly supported functions in $\R^d$ (i.e not as periodic functions). Then we can construct a family of \emph{concentrated Mikado densities and fields}, parametrized by a concentration parameter $\mu > 0$ , defined as a rescaled version of $\Theta$ and $W$, as follows:
\begin{equation*}
\Theta_\mu (x) := \mu^\alpha \Theta(\mu x), \qquad W_\mu(x) = \mu^\beta W(\mu x).
\end{equation*} 
It is now not difficult to see that, if \eqref{eq:condition} is satisfied, then one can choose $\alpha$, $\beta$ so that 
\begin{equation}
\label{eq:order1}
\|\Theta_\mu\|_{L^p} \approx 1, \qquad \|W_\mu\|_{L^{p'}} \approx 1,
\end{equation}
and
\begin{equation}
\label{eq:derivative:small}
\|D W_\mu\|_{L^{\tilde p}} \approx \mu^{-c},
\end{equation}
for some  $c > 0$, so that $\|D W_\mu\|_{L^{\tilde p}} \to 0$, as $\mu \to \infty$. In this way, we can produce a whole family of Mikado fields, which ``are not degenerating'' as $\mu \to \infty$ (i.e. they remains ``of order $1$'', in some suitable norm, thanks to \eqref{eq:order1}), but, at the very same time, have vanishing derivative, thanks to \eqref{eq:derivative:small}. 

We can now modify our \emph{Ansatz} \eqref{eq:induction} as follows:
\begin{equation}
\label{eq:new:induction}
\rho_{q+1} = \rho_q + a_q(t,x) \Theta_{\mu_q}(\lambda_q x), \qquad u_{q+1} = u_q + b_q(t,x) W_{\mu_q}(\lambda_q x),
\end{equation}
where $\mu_q$ is a sequence of real numbers, with $\mu_q \to \infty$ as $q \to \infty$, to be chosen appropriately. In this way, thanks to \eqref{eq:derivative:small}, the estimate in \eqref{eq:estimate:der} becomes
\begin{equation*}
\|D u_{q+1} - Du_q\|_{L^{\tilde p}} \approx \lambda_q \|b_q\|_{L^\infty} \|DW_{\mu_q}\|_{L^{\tilde p}} \lesssim  \|b_q\|_{L^\infty} \lambda_q \mu_q^{-c},
\end{equation*} 
and thus, if $\mu_q$ is chosen much bigger than $\lambda_q$, the distance $\|D u_{q+1} - Du_q\|_{L^{\tilde p}} $ can be made arbitrarily small, thus getting convergence of $u_q$ in $W^{1, \tilde p}$ and hence proving Theorem \ref{thm:main}.

% You may incorporate your references as follows in your main tex file.
% Using BibTex is not recommended but can be handled.

\end{document}